\documentclass{article}%
\usepackage{amssymb}
\usepackage{graphicx}
\usepackage{amsmath}
\usepackage{amsfonts}%
\providecommand{\U}[1]{\protect\rule{.1in}{.1in}}
\newtheorem{theorem}{Theorem}

\newtheorem{corollary}[theorem]{Corollary}

\newtheorem{lemma}[theorem]{Lemma}

\newtheorem{proposition}[theorem]{Proposition}
\newtheorem{remark}[theorem]{Remark}

\begin{document}

\begin{center}
{\large On the local limit theorems for psi-mixing Markov chains.}

\bigskip

Florence Merlev\`{e}de,\ Magda Peligrad and Costel Peligrad

\bigskip

\end{center}

Universit\'{e} Gustave Eiffel, LAMA and CNRS UMR 8050. 

Email: florence.merlevede@univ-eiffel.fr

Department of Mathematical Sciences, University of Cincinnati, PO Box 210025,
Cincinnati, Oh 45221-0025, USA. \texttt{ }

Email: peligrm@ucmail.uc.edu, peligrc@ucmail.uc.edu

\bigskip

\noindent\textbf{Keywords}: \noindent Markov chains, local limit theorem, mixing.

\bigskip

\noindent\textbf{2010 Mathematics Subject Classification}: Primary 60F05; 60J05

\begin{center}
\bigskip\bigskip

Abstract
\end{center}

In this paper we investigate the local limit theorem for additive functionals
of a nonstationary Markov chain with finite or infinite second moment. The
moment conditions are imposed on the individual summands and the weak
dependence structure is expressed in terms of some uniformly mixing coefficients.

\section{Introduction}

A local limit theorem for partial sums $(S_{n})_{n\geq1}$ of a sequence of
centered random variables is a result about the rate of convergence of the
probabilities of the type $P(a\leq S_{n}\leq b)$. Local limit theorems have
been studied for\ the case of lattice random variables and the case of
non-lattice random variables. The lattice case means that there exists $v>0$
and $a\in\mathbb{R}$ such that the values of all the variables in the sum
$S_{n}$ are concentrated on the lattice $\{a+kv:k\in\mathbb{Z}\}$, whereas the
nonlattice case means that no such $a$ and $v$ exists. This problem was
intensively studied for i.i.d. sequences of random variables starting with
works by Shepp (1964), Stone (1965) and Feller (1967). When the variables are
independent but not necessarily identically distributed we mention papers by
Mineka and Silverman (1970), Maller (1978), Shore (1978). Dolgopyat (2016)
treated the vector valued sequences of independent random variables. In the
lattice case, for countable state Markov chains with finite second moments,
the local limit theorem is discussed in Nagaev (1963) and S\'{e}va (1995)
while the case of infinite variance is analyzed in Aaronson and Denker
(2001a,b) and Szewczak (2010).

In the lattice case, Szewczak (2010) established a local limit theorem for
functions of continued fractions expansions, which is an example of $\psi
-$mixing sequence. Also in the stationary case we mention the local limit
\ theorems for Markov chains in the papers by Herv\'{e} and P\`{e}ne (2010),
Ferr\'{e} et al. (2012). Hafouta and Kifer (2016) proved a local limit theorem
for nonconventional sums of stationary a class of uniformly mixing Markov
chains satisfying conditions related to $\psi-$mixing chains. As noticed in
Denker (1992), Bryc (1990,1992) the concept of $\psi-$mixing is well suited to
derive large deviation results. As examples of $\psi-$mixing Markov chains we
mention Gibbs-Markov dynamical systems introduced in Aaronson and Denker
(2001b), which contain finite state aperiodic Markov chains and certain
recurrent Markov chains with infinite state space.

We assume now that $(\xi_{k})_{k\geq1}$ is a Markov chain defined on
$(\Omega,\mathcal{K},P)$ with values in $(S,\mathcal{B}(S))$ with regular
transition probabilities,
\[
Q_{k}(x,A)=P(\xi_{k}\in A|\xi_{k-1}=x);\text{ }P_{k}(A)=P(\xi_{k}\in A).
\]
Also denote by $Q_{k}$ the associated operator defined on $L_{2}%
(\mathcal{B}(S))$ by $Q_{k}g(x)=\int g(y)Q_{k}(x,dy)$.

For some real-valued measurable functions $g_{j}$ on $S$ define
\begin{equation}
X_{j}=g_{j}(\xi_{j}). \label{defX}%
\end{equation}
The sequence $(X_{j})_{j\geq1}$ is assumed centered $(E(X_{j})=0$ for all
$j\geq1)$ and, unless otherwise specified, having finite second moments
$(E(X_{j}^{2})<\infty$ for all $j\geq1)$. Set%
\[
S_{n}=\sum\nolimits_{k=1}^{n}X_{k}\text{, }\sigma_{n}^{2}=E(S_{n}^{2})\text{
and }\tau_{n}^{2}=\sum_{j=1}^{n}E(X_{j}^{2}).
\]
We shall assume that there are two constant $a>0$ and $b<\infty$ with the
following property:

For all $k\geq2$ there is $S_{k}^{\prime}\in\mathcal{B}(S)$ with
$P_{k-1}(S_{k}^{\prime})=1$ such that for all $A\mathcal{\in B}(S)$ and $x\in
S_{k}^{\prime}$ we have
\begin{equation}
aP_{k}(A)\leq Q_{k}(x,A)\leq bP_{k}(A)\text{.}\label{cond0}%
\end{equation}

Denote%
\[
\gamma=\frac{a^{4}}{b}\text{ }.
\]
Clearly $b\geq1$ and $a\leq1.$

Throughout the paper we shall assume that $\tau_{n}^{2}\rightarrow\infty.$ As
we shall see latter, since we assume $a>0,$ the condition $\tau_{n}%
^{2}\rightarrow\infty$ is equivalent to $\sigma_{n}^{2}\rightarrow\infty$ (see
subsection \ref{subsectionmixingvar}).

In order to extend the local CLT beyond the nonstationary sequences of i.i.d.,
we shall combine several techniques specifically designed for obtaining local
limit theorems with a bound on the characteristic function using $\gamma.$
More precisely, the conditions and techniques are rooted in Mineka and
Silverman (1970) and Maller (1978), who treated the local limit theorem for
nonstationary sequences of i.i.d. in the non-lattice case. We shall prove that
if we assume (\ref{cond0}),$~$then the results referring to the local CLT in
Maller (1978) and also Mineka and Silverman (1970) can be extended from
independent sequences of random variables to the Markovian case. Furthermore
we shall also consider the situation when the variables have infinite variance
and are in the domain of attraction of the normal law.

Our paper is organized as follows. In Section 2 we present the local limit
theorem for nonstationary Markov chains. In section Proofs we present bounds
on the characteristic function of sums, bounds of the variance of sums and the
proof of the main results.

In the following section, the notation $a(n)=o(n)$ means that
$a(n)/n\rightarrow0$ as $n\rightarrow\infty.$ Also by $\Rightarrow$ we denote
the convergence in distribution.

\section{Local central limit theorem}

In the sequel we shall denote by $f_{k}(t)$ the Fourier transform of $X_{k},$%
\[
f_{k}(t)=f_{X_{k}}(t)=E(\exp(itX_{k})).
\]

\bigskip

The first condition is the usual Lindeberg condition and is imposed in order
to obtain the CLT.

\bigskip

\textbf{Lindeberg's condition.} For any $\varepsilon>0,$
\begin{equation}
\frac{1}{\tau_{n}^{2}}\sum\limits_{k=1}^{n}E(X_{k}^{2}I(|X_{k}|\geq
\varepsilon\tau_{n}))\rightarrow0. \label{lindeberg1}%
\end{equation}
The second condition deals with the characteristic functions of normalized
individual summands on bounded but large intervals.

\bigskip

\textbf{Condition A}. There is $\delta>0$ and $n_{0}\in N$ such that for
$1\leq|u|\leq\delta\tau_{n}$ and $n>n_{0}$
\begin{equation}
\frac{\gamma}{8}\sum\limits_{k=1}^{n}(1- \big |f_{k}(\frac{u}{\tau_{n}})
\big |^{2})>g(u)\text{ and }\exp(-g(u))\text{ is integrable on }%
\mathbb{R}\text{.} \label{regular}%
\end{equation}
Next condition is a nonlattice-type condition. Recall that a random variable
$X$ has a nonlattice distribution is equivalent to $|f_{X}(t)|<1$ for all $t$.

\bigskip

\textbf{Condition B}. For $u\neq0$ there is $c(u)$, an open interval $O_{u}$
containing $u$ and a $n_{0}=n_{0}(u)$ such that for all $t\in O_{u}$, and
$n>n_{0}$
\begin{equation}
\frac{\gamma}{8(\ln\tau_{n})}\sum\limits_{k=1}^{n}(1-|f_{k}(t)|^{2})\geq
c(u)>1. \label{verC2}%
\end{equation}
Our general local limit theorem is as follows:

\begin{theorem}
\label{ThLocal-lattice}Let $(X_{j})_{j\geq1}$ be defined by (\ref{defX}).
Assume that Conditions A , B, (\ref{cond0})\ and (\ref{lindeberg1}) are
satisfied. Then, for any function $h$ on $\mathbb{R}$ which is continuous and
with compact support,%
\begin{equation}
\lim_{n\rightarrow\infty}\sup_{u\in\mathbb{R}} \big |\sqrt{2\pi}\sigma
_{n}Eh(S_{n}-u)-\exp(-u^{2}/2\sigma_{n}^{2})\int h(u)\lambda(du) \big |=0,
\label{LCLTnonlat}%
\end{equation}
where $\lambda$ is the Lebesgue measure.
\end{theorem}

It is well known that the convergence in (\ref{LCLTnonlat}) implies that for
any $c$ and $d$ real numbers with $c<d$
\[
\lim_{n\rightarrow\infty}\sup_{u\in\mathbb{R}}\bigg|\sqrt{2\pi}\sigma
_{n}P(c+u\leq S_{n}\leq d+u)-(d-c)\exp(-u^{2}/2\sigma_{n}^{2})\ \bigg|=0.
\]
In particular, since $\sigma_{n}\rightarrow\infty$ as $n\rightarrow\infty$,
then for fixed $A>0$,
\[
\lim_{n\rightarrow\infty}\sup_{|u|\leq A}\bigg|\sqrt{2\pi}\sigma_{n}P(c+u\leq
S_{n}\leq d+u)-(d-c)\bigg|=0.
\]
If we further take $u=0$, then
\[
\lim_{n\rightarrow\infty}\sqrt{2\pi}\sigma_{n}P(S_{n}\in\lbrack c,d])=d-c.
\]
In other words, the sequence of measures $\sqrt{2\pi}\sigma_{n}P(S_{n}%
\in\lbrack c,d])$ of the interval $[c,d]$ converges to Lebesgue measure.

\begin{remark}
Theorem \ref{ThLocal-lattice} can be reformulated for triangular arrays of
Markov chains $(\xi_{n,i})_{1\leq i\leq n}$ and $X_{n,i}=g_{n,i}(\xi_{n,i})$.
The difference is that in condition (\ref{cond0}) and in Conditions A and B
the quantities will depend on $n$. Namely, now we have $a_{n},$ $b_{n}%
,\gamma_{n},$ and $f_{n,k}.$ Also the Lindeberg's condition should be adjusted
as in Theorem 1 of Gudynas (1977) which is a variant of the Dobrushin's (1956)
CLT\ for nonstationary Markov chains. So we have to replace (\ref{lindeberg1})
by
\[
\frac{1}{a_{n}\tau_{n}^{2}}\sum\limits_{k=1}^{n}E(X_{n,k}^{2}I(|X_{n,k}%
|\geq\varepsilon a_{n}\tau_{n}))\rightarrow0\text{ as }n\rightarrow\infty.
\]

\end{remark}

\textbf{Discussion on Conditions A and B:}

\bigskip

By using Condition (1.1) in Maller (1978) and Lindeberg's condition, Condition
A can be verified under the following "balance" type condition (its proof is
postponed to the end of the paper in Lemma \ref{LemmaA}):

\bigskip

\textbf{Condition A}$_{1}$. There is $0\leq c<1$ and $\delta>0$ such that%
\[
\lim\sup_{n\rightarrow\infty}\frac{\sum\nolimits_{k=1}^{n}E(X_{k}^{2}%
I(|X_{k}|>\delta))}{\tau_{n}^{2}}<c.
\]
As we can easily verify, condition A$_{1}$ is satisfied under stronger
condition: There is $0\leq c<1$ and $\delta>0$ such that%
\begin{equation}
\frac{E(X_{k}^{2}I(|X_{k}|>\delta))}{E(X_{k}^{2})}<c\text{ for all }k,
\label{to verify cond A}%
\end{equation}
and also under Mineka and Silverman (1970) condition, namely: For some
$\delta>0$ and $d>0$ ,
\[
E(X_{k}^{2}I(|X_{k}|\leq\delta))\geq dE(X_{k}^{2})\text{ for all }k.
\]
Clearly Condition A$_{1}$ is trivially implied if there is $C>0$ such that
$|X_{k}|\leq C$ a.s. or under the near stationarity assumption:

There is a random variable $X$ and constants $0<c_{1}\leq1$ and $c_{2}\geq1$
such that%
\[
c_{1}P(|X|\geq x)\leq P(|X_{k}|\geq x)\leq c_{2}P(|X|\geq x)\text{ for all
}x\geq1\text{ and all }k\in\mathbb{N}.
\]

\bigskip

Condition B is satisfied under condition (1.3) of Mineka and Silverman (1970)
(which is an adaptation of a condition due to Rozanov, 1957):

\bigskip

\textbf{Condition B}$_{1}$. For $u\neq0$ there is an $\varepsilon
=\varepsilon(u)>0,$ for which%
\[
\frac{1}{\ln\tau_{n}}\sum\limits_{j=1}^{n}P(X_{j}-a_{j}\in A(u,\varepsilon
))\rightarrow\infty,
\]
where $a_{i}$ is a bounded sequence of constants satisfying $\inf_{1\leq
j\leq\infty}P(|X_{j}-a_{j}|<\delta)>0$ for every $\delta>0$ and
$A(u,\varepsilon)=\{x:|x|<M,|xu-\pi m|\geq\varepsilon\},$ for each integer $m$
with $|m|\leq M,$ where $M>0$ is fixed, large enough such that $\inf_{1\leq
j\leq\infty}P(|X_{j}|<M)>0$ (the existence for such an $M$ is a part of the assumption).

The fact that Condition B$_{1}$ implies Condition B was proven by Mineka and
Silverman (1970). Under Condition B$_{1}$, Mineka and Silverman (1970), on the
top of page 595, showed that for each $u$ there is a positive constant $K_{u}$
independent on $k$ such that for all $t$ such that $|t-u|<\varepsilon/4M,$ and
for all $k\in\mathbb{N}$ we have
\[
|f_{k}(t)|^{2}-1\leq-\frac{1}{4}K_{u}\varepsilon^{2}P(X_{k}\in A(u,\varepsilon
)).
\]
Also, from Corollary 1 in Mineka and Silverman (1970), Condition B$_{1}$ can
be replaced with the stronger condition: the variables $X_{k}$'s have
uniformly bounded densities, or by Corollary 2 in the same paper, Condition
B$_{1}$ is satisfied if there are three rationally independent numbers
$d_{1},d_{2},d_{3}$ such that $\inf_{1\leq k\leq\infty}P(|X_{k}-d_{j}%
|<\delta)>0$ for $j=1,2,3.$

\bigskip

By using stronger degrees of stationarity, we can further simplify Condition B:

\textbf{Condition B}$_{2}.$ Assume (\ref{lindeberg1})\ and for all $u\neq0$
there is an open interval $O_{u}$ containing $u$ and a $n_{0}=n_{0}(u)$ such
that for all $t\in O_{u}$ and $n>n_{0},$
\begin{equation}
\frac{1}{n}\sum\limits_{k=1}^{n}|f_{k}(t)|^{2}<1. \label{to verify B 2}%
\end{equation}
In order to see that Condition B$_{2}$ implies Condition B, we note that by
Lindeberg's condition (\ref{lindeberg1})\ we obtain
\[
\lim_{n\rightarrow\infty}\frac{\tau_{n+1}}{\tau_{n}}=1.
\]
This gives that, for any $d>1$ and all $n$ sufficiently large
\[
\frac{\tau_{n}}{\tau_{n_{0}}}=\frac{\tau_{n}}{\tau_{n-1}}\cdot\frac{\tau
_{n-1}}{\tau_{n-2}}\cdot...\cdot\frac{\tau_{n_{0}+1}}{\tau_{n_{0}}}\leq
d^{n-n_{0}}.
\]
Therefore,
\[
\ln\tau_{n}\leq(n-n_{0})\ln d+\ln\tau_{n_{0}}=n\ln d+o(1).\text{ }%
\]
and then%
\[
\frac{1}{\ln\tau_{n}}\sum\limits_{k=1}^{n}(1-|f_{k}(t)|^{2})\geq\frac{1}{\ln
d+o(1)}(1-\frac{1}{n}\sum\limits_{k=1}^{n}|f_{k}(t)|^{2}).
\]
So, (\ref{verC2}) is satisfied if we show that for $|t-u|\leq\varepsilon,$ we
have:
\[
\frac{1}{\ln d+o(1)}(1-\frac{1}{n}\sum\limits_{k=1}^{n}|f_{k}(t)|^{2}%
)\geq\frac{8}{\gamma}.
\]
This is equivalent to showing that for $t$ such that $|t-u|\leq\varepsilon$
and $n>n_{0}$ we have%
\[
\frac{1}{n}\sum\limits_{k=1}^{n}|f_{k}(t)|^{2}\leq1-\frac{8}{\gamma}\ln
d+o(1).
\]
If we select now $d>1$ close enough to $1$, we see that Condition B$_{2}$
implies Condition B.

Condition B$_{2}$ is satisfied in the stationary case if the marginal
distribution satisfies $|f_{0}(t)|<1$ for all $t.$ It is well-known that
$|f_{0}(t)|<1$ for all $t\neq0$ is equivalent to $X_{0}$ not having a lattice distribution.

By using Theorem \ref{ThLocal-lattice} we can treat linear statistics with
coefficients which are uniformly bounded above and away from zero:

\begin{corollary}
\label{corsta-lin}Assume that $(\xi_{k})_{k\in\mathbb{Z}}$ is a strictly
stationary Markov chain. For a measurable function $g$ and $k\in\mathbb{Z}$
define $X_{k}=g(\xi_{k})$ Assume that $E(X_{k})=0$ and $\ E(X_{k}^{2}%
)<\infty.$ Assume that $X_{0}$ has a non-lattice distribution and condition
(\ref{cond0}) is satisfied. We consider an array of real numbers
$(a_{n,k})_{k\geq1}$ such that there are two positive constants $m,M$ with
$0<m\leq|a_{n,k}|\leq M$ for all $n$ and $k.$ Define%
\[
X_{n,\ell}=a_{n,\ell}X_{\ell}.
\]
Then, for any function $h$ on $\mathbb{R}$ which is continuous and with
compact support, $S_{n}=\sum\nolimits_{k=1}^{n}X_{n,\ell}$ satisfies Theorem
\ref{LCLTnonlat}.
\end{corollary}

\begin{remark}
Note that the strictly stationary case follows from Corollary \ref{corsta-lin}
if we take for all $1\leq\ell\leq n$ the constants $a_{n,\ell}=1.$
\end{remark}

With a very similar proof as of Corollary \ref{corsta-lin} we can treat the
linear processes with short memory.

\begin{corollary}
\label{corlinshort}Let $(X_{k})_{k\in\mathbb{Z}}$ be as in Corollary
\ref{corsta-lin}. Let $(a_{i})_{i\geq1}$ be a sequence of real numbers such
that $\sum_{i\geq1}|a_{i}|<\infty$. and $m=\inf_{j}|A_{j}|>0,\ $where
$A_{j}=a_{1}+a_{2}+...+a_{j}.$ Construct%
\[
Y_{k}=\sum_{i\geq1}a_{i}X_{k+i}\text{ and }S_{n}=\sum\limits_{k=1}^{n}Y_{k}.
\]
Assume that $X_{0}$ has a nonlattice distribution and condition (\ref{cond0})
is satisfied. Then
\[
\lim_{n\rightarrow\infty}\sup_{u\in\mathbb{R}} \big |\sqrt{2\pi}\sigma_{n}|A
|Eh(S_{n}-u)-\exp(-u^{2}/2\sigma_{n}^{2}A^{2})\int h(u)\lambda(du) \big |=0,
\]
where $v_{n}^{2}=E(\sum_{i=1}^{n}X_{i})^{2}$ and $A=\lim_{j\rightarrow\infty
}A_{j}.$
\end{corollary}

\bigskip

We can also provide a result for the stationary situation when the variance of
the individual summands can be infinite. As an application of the proof of
Theorem \ref{ThLocal-lattice} we obtain the following corollary:

\begin{corollary}
\label{corstationary infvar}Assume that $(\xi_{k})_{k\in\mathbb{Z}}$ is a
strictly stationary Markov chain. Define $(X_{k})_{k\in\mathbb{Z}}$ by
$X_{k}=g(\xi_{k})$ and assume $E(X_{0})=0$ and $H(x)=E(X_{0}^{2}I(|X_{0}|\leq
x)$ is a slowly varying function as $x\rightarrow\infty$. Assume (\ref{cond0})
and $X_{0}$ has a non-lattice distribution. Then there is $b_{n}%
\rightarrow\infty$ such that for any function $h$ on $\mathbb{R}$ which is
continuous and with compact support,
\[
\lim_{n\rightarrow\infty}\sup_{u\in\mathbb{R}}|\sqrt{2\pi}b_{n}Eh(S_{n}%
-u)-\exp(-u^{2}/2b_{n}^{2})\int h(u)\lambda(du)|=0.
\]

\end{corollary}

As far as we know this corollary is new, though for Gibbs-Markov processes the
result is contained in Aaronson and Denker (2001a) and for continued fraction
processes can be found in Szewczak (2010).

\bigskip

\textbf{Example.} For every irrational number $x$ in $(0,1)$ there is a unique
sequence of positive integers $x_{1},x_{2},x_{3},...$ such that the following
continued fraction expansion holds:%
\[
x=\frac{1}{x_{1}+\frac{1}{x_{2}+\frac{1}{x_{3+\cdot\cdot\cdot}}}}.
\]
If we introduce on $[0,1]$ the Gauss probability measure with the density
$f(x)=(\ln2)^{-1}(1+x)^{-1},$ then the sequence $(x_{1},x_{2},x_{3},...)$ is a
strictly stationary Markov chain. We know from Lemma 2.1 in Philipp (1988)
that one can take for $a$ and $b$ which appear in condition \ref{cond0}
$a=0.2$ and $b=1.8.$ Now we can consider $X_{k}$ defined by (\ref{defX}) as a
measurable function of $x_{k}$ and apply both Corollary \ref{corsta-lin} and
Corollary \ref{corstationary infvar}.

\section{Proofs}

\subsection{Bounds on the characteristic function}

The bound on the characteristic function of a Markov chain is inspired by
Lemma 1.5 in Nagaev (1961). It is given in the following proposition.

\begin{proposition}
\label{factor}Let $(X_{j})_{j\geq1}$ be defined by (\ref{defX}). Then%
\[
|E(\exp(iuS_{n})|^{4}\leq\prod\nolimits_{j=1}^{n}[1-\frac{\gamma}{2}%
(1-|f_{j}(u)|^{2})]\text{ .}%
\]

\end{proposition}

For proving this proposition we need some preliminary considerations. For $u$
fixed let us introduce the operator $T_{k}=T_{u,k}$ defined on complex valued
bounded functions by:
\[
T_{k}(h)(x)=\int h(y)\exp(iug_{k}(y))Q_{k}(x,dy).
\]
So%
\[
T_{k}(h)(\xi_{k-1})=E([h(\xi_{k})\exp(iuX_{k})]|\xi_{k-1}).
\]
Notice that the values are also complex bounded functions.

For an operator $T$ on $\mathbb{L}_{\infty}(S,\mathcal{B}(S))$ denote by
$||T||$=$\sup_{|f|_{\infty}<1}|T_{k}(f)|_{\infty}$.

\begin{lemma}
\label{estimate3}For any $k\in N,$ $u\in\mathbb{R}$ we have for all $k\geq2,$
\[
||T_{k-1}\circ T_{k}||^{2}\leq1-\frac{\gamma}{2}(1-|E(\exp(iuX_{k-1}%
)|^{2})\text{ .}%
\]

\end{lemma}

\noindent\textbf{Proof.} Let $x\in S^{\prime},$ where $S^{\prime}\in$
$\mathcal{B}(S)$ such that $P_{k-1}(S^{\prime})=1,$ for which condition
(\ref{cond0}) holds. By the definition of $T_{k}$'s
\begin{gather*}
T_{k-1}\circ T_{k}(h)(x)=\\
=\int\exp(iug_{k-1}(y))\int h(z)\exp(iug_{k}(z))Q_{k}(y,dz)Q_{k-1}(x,dy).
\end{gather*}
Changing the order of integration
\begin{align*}
T_{k-1}\circ T_{k}(h)(x)  &  =\int h(z)\exp(iug_{k}(z))\int\exp(iug_{k-1}%
(y))Q_{k-1}(x,dy)Q_{k}(y,dz)\\
&  =\int h(z)\exp(iug_{k}(z))m_{x}(dz)
\end{align*}
where, for $x$ fixed $m_{x}$ is measure defined on $\mathcal{B}(S)$ by
\[
m_{x}(A)=\int\exp(iug_{k-1}(y))Q_{k}(y,A)Q_{k-1}(x,dy).
\]
Denote by $\mathrm{Var}\left(  m_{x}\right)  $ the total variation of $m_{x}.$
With this notations and because $h$ is bounded by $1$,%
\[
|T_{k-1}\circ T_{k}(h)(x)|\leq\mathrm{Var}\left(  m_{x}\right)  .
\]
Now, in order to compute the total variation for $m_{x}$ we start from the
following estimate%
\begin{gather*}
\left(  \int Q_{k}(y,A)Q_{k-1}(x,dy)\right) ^{2}- \Big |\int\exp
(iug_{k-1}(y))Q_{k}(y,A)Q_{k-1}(x,dy) \Big |^{2}=\\
\iint(1-\cos(u\left(  g_{k-1}(y)-g_{k-1}(y^{\prime}))\right)  Q_{k}%
(y,A)Q_{k-1}(x,dy)Q_{k}(y^{\prime},A)Q_{k-1}(x,dy^{\prime})\\
\iint2\sin^{2}\left(  \frac{u}{2}(g_{k-1}(y)-g_{k-1}(y^{\prime}))\right)
Q_{k}(y,A)Q_{k-1}(x,dy)Q_{k}(y^{\prime},A)Q_{k-1}(x,dy^{\prime})\\
\geq a^{4}P_{k}^{2}(A)\iint2\sin^{2}\left(  \frac{u}{2}(g_{k-1}(y)-g_{k-1}%
(y^{\prime}))\right)  P_{k-1}(dy)P_{k-1}(dy^{\prime}).
\end{gather*}
But%
\[
\int Q_{k}(y,A)Q_{k-1}(x,dy)+ \Big |\int\exp(iug_{k-1}(y))Q_{k}(y,A)Q_{k-1}%
(x,dy) \Big |\leq2bP_{k}(A).
\]
So, by dividing the last two inequalities, we obtain
\begin{gather*}
\left(  \int Q_{k}(y,A)Q_{k-1}(x,dy)\right)  - \Big |\int\exp(iug_{k-1}%
(y))Q_{k}(y,A)Q_{k-1}(x,dy) \Big |\geq\\
\frac{a^{4}}{2b}P_{k}(A)\iint2\sin^{2}\left(  \frac{u}{2}(g_{k-1}%
(y)-g_{k-1}(y^{\prime}))\right)  P_{k-1}(dy)P_{k-1}(dy^{\prime})\\
=\frac{a^{4}}{2b}P_{k}(A)(1-|f_{k-1}(u)|^{2}).
\end{gather*}
So%
\[
|m_{x}(A)|\leq\int Q_{k}(y,A)Q_{k-1}(x,dy)-\frac{\gamma}{2}P_{k}%
(A)(1-|f_{k-1}(u)|^{2}).
\]
Now we consider $(A_{i})_{i\in J}$ a finite partition\ of $S,$ with sets in
$\mathcal{B}(S)$. Then%
\[
\sum\nolimits_{i\in J}|m_{x}(A_{i})|\leq1-\frac{\gamma}{2}(1-|f_{k-1}%
(u)|^{2}).
\]
It follows that, for all $x\in S^{\prime}$
\[
\mathrm{Var}\left(  m_{x}\right)  \leq1-\frac{\gamma}{2}(1-|f_{k-1}%
(u)|^{2})\text{ }.
\]
and Lemma \ref{estimate3} follows. $\square$

\bigskip

\noindent\textbf{Proof of Proposition \ref{factor}.} Note that
\[
E(\exp(iuS_{2k})|\xi_{0}=x)=T_{1}\circ T_{2}\circ...\circ T_{2k}(1)(x).
\]
So
\[
|E(\exp(iuS_{2k})|\xi_{0})|\leq||T_{1}\circ T_{2}||...||T_{2k-1}\circ
T_{2k}||\text{ \ a.s.}%
\]
By Lemma \ref{estimate3} we have that, for $k\geq1,$%
\[
|E(\exp(iuS_{2k})|\xi_{0})|^{2}\leq\prod\nolimits_{j=1}^{k}[1-\frac{\gamma}%
{2}(1-|f_{2j-1}(u)|^{2})]\text{ a.s.}%
\]
Also, by Lemma \ref{estimate3}, for $k\geq1,$
\begin{gather*}
|E(\exp(iuS_{2k})|\xi_{0})|^{2}\leq||T_{2}\circ T_{3}||^{2}...||T_{2k-2}\circ
T_{2k-1}||^{2}||T_{2k}(1)||^{2}\\
\leq\prod\nolimits_{j=1}^{k}[1-\frac{\gamma}{2}(1-|f_{2j}(u)|^{2})]\text{
\ a.s. }%
\end{gather*}
and so, by multiplying these two relations we get%
\[
|E(\exp(iuS_{2k})|\xi_{0})|^{4}\leq\prod\nolimits_{j=1}^{2k}[1-\frac{\gamma
}{2}(1-|f_{j}(u)|^{2})]\text{ a.s.}%
\]
A similar result can be obtain for $|E(\exp(iuS_{2k+1})|\xi_{0})|^{4}$. The
result in Proposition \ref{factor} follows. $\square$

\subsection{Mixing conditions and the variance of partial sums}

\label{subsectionmixingvar}

We shall clarify here the relation between $a$ and $b$ in condition
(\ref{cond0}) and several mixing coefficients for stochastic processes. Let
$(\Omega,\mathcal{K},P)$ be a probability space and let $\mathcal{A}%
,\mathcal{B}$ be two sub $\sigma$-algebras of $\mathcal{K}$. Define the
maximal coefficient of correlation
\[
\rho(\mathcal{A},\mathcal{B})=\sup_{X\in L_{2}(\mathcal{A}),Y\in
L_{2}(\mathcal{B})}|\mathrm{corr}\,(X,Y)|\text{ ,}%
\]
where $L_{2}(\mathcal{A})$ is the space of random variables that are
$\mathcal{A}$ measurable and square integrable.

Relevant to our paper are the lower and upper $\psi-$mixing coefficients
defined by
\begin{align*}
\psi^{\prime}(\mathcal{A},\mathcal{B})  &  =\inf\frac{P(A\cap B)}%
{P(A)P(B)};\text{ }A\in\mathcal{A}\text{ and }B\in\mathcal{B}\text{,
}P(A)P(B)>0.\\
\psi^{\ast}(\mathcal{A},\mathcal{B})  &  =\sup\frac{P(A\cap B)}{P(A)P(B)}%
;\text{ }A\in\mathcal{A}\text{ and }B\in\mathcal{B}\text{, }P(A)P(B)>0.
\end{align*}
By a result of Bradley (2020) we have the following lemma:

\begin{lemma}
\label{Bradley2020}(Bradley, 2020)%
\begin{equation}
\rho(\mathcal{A},\mathcal{B})\leq1-\psi^{\prime}(\mathcal{A},\mathcal{B}).
\label{rel ro}%
\end{equation}

\end{lemma}

Proof. For simplicity we denote $\rho=\rho(\mathcal{A},\mathcal{B})$ and
$\psi^{\prime}=\psi^{\prime}(\mathcal{A},\mathcal{B}).$ Without restricting
the generality we assume $\psi^{\prime}>0.$ By the definition of $\rho$ we
have to show that for any $X\in L_{2}(\mathcal{A})$ and $Y\in L_{2}%
(\mathcal{B})$, we have to show that
\[
|E(XY)|\leq(1-\psi^{\prime})||X||_{2}||Y||_{2}.
\]
By a measure theoretic argument, for variables with values in a separable
Hilbert space, it is enough to prove this lemma for simple functions with mean
zero. So, let $X=\sum\nolimits_{i=1}^{n}a_{i}I(A_{i})$ and $Y=\sum
\nolimits_{j=1}^{m}b_{j}I(B_{j}),$ where $A_{i}\in\mathcal{A}$ and $B_{j}%
\in\mathcal{B}$ are partitions of $\Omega$ and $X$ and $Y$ have mean $0.$

Denote
\[
R(A_{i}B_{j})=(1-\psi^{\prime})^{-1}P(A_{i}B_{j})-(1-\psi^{\prime})^{-1}%
\psi^{\prime}P(A_{i})P(B_{j})
\]
and note that, by the definition of $\psi^{\prime},$ for all $i$ and $j$ we
have that $R(A_{i}B_{j})\geq0.$ Also \
\[
\sum\nolimits_{i=1}^{n}R(A_{i}B_{j})=P(B_{j})\text{ and }\sum\nolimits_{j=1}%
^{m}R(A_{i}B_{j})=P(A_{i}).
\]
Moreover, we have the decomposition%
\[
P(A_{i}B_{j})=\psi^{\prime}P(A_{i})P(B_{j})+(1-\psi^{\prime})R(A_{i}B_{j}).
\]
Now, since $E(X)=0,$ clearly $\sum\nolimits_{i=1}^{n}a_{i}P(A_{i})=0$ and
therefore, by the above identity,%
\begin{align*}
E(XY)  &  =\sum\nolimits_{i,j}a_{i}b_{j}P(A_{i}B_{j})=\sum\nolimits_{i,j}%
a_{i}b_{j}\left(  \psi^{\prime}P(A_{i})P(B_{j})+(1-\psi^{\prime})R(A_{i}%
B_{j})\right) \\
&  =(1-\psi^{\prime})\sum\nolimits_{i,j}a_{i}b_{j}R(A_{i}B_{j}).
\end{align*}
It follows that%
\[
|E(XY)|\leq(1-\psi^{\prime})\sum\nolimits_{i,j}|a_{i}b_{j}|R(A_{i}B_{j}).
\]
So, by applying Holder's inequality twice,%
\begin{align*}
\sum\nolimits_{i,j}|a_{i}b_{j}|R(A_{i}B_{j})  &  \leq\sum\nolimits_{i}%
|a_{i}|\ \left(  \sum\nolimits_{j}R(A_{i}B_{j})\right)  ^{1/2}\left(
\sum\nolimits_{j}b_{j}^{2}R(A_{i}B_{j})\right)  ^{1/2}\\
\leq &  \left[  \sum\nolimits_{i}\sum\nolimits_{j}|a_{i}|^{2}R(A_{i}%
B_{j})\right]  ^{1/2}\left[  \sum\nolimits_{i}\sum\nolimits_{j}|b_{j}%
|^{2}R(A_{i}B_{j})\right]  ^{1/2}\\
&  =\left[  \sum\nolimits_{i}|a_{i}|^{2}P(A_{i})\right]  ^{1/2}\left[
\sum\nolimits_{j}|b_{j}|^{2}P(B_{j})\right]  ^{1/2}=||X||_{2}||Y||_{2}.
\end{align*}
$\ \square$

\bigskip

For a sequence $\mathbf{X}=(X_{k})_{k\geq1}$ of random variables $\psi
_{k}^{\prime}(\mathbf{X)}=\inf_{m\geq1}\psi^{\prime}(\mathcal{F}_{1}%
^{m},\mathcal{F}_{k+m}^{\infty})$, $\psi_{k}^{\ast}(\mathbf{X)}=\sup_{m\geq
1}\psi^{\ast}(\mathcal{F}_{1}^{m},\mathcal{F}_{k+m}^{\infty})$ and $\rho
_{k}(\mathbf{X)=}\sup_{m\geq1}\rho(\mathcal{F}_{1}^{m},\mathcal{F}%
_{k+m}^{\infty}),$ where $\mathcal{F}_{k}^{m}=\sigma(X_{j},k\leq j\leq m)$.

For a Markov chain $\mathbf{\xi}=(\xi_{k})_{k\geq1}$ the definitions simplify
\begin{align*}
\psi_{k}^{\prime}  &  =\inf_{m\geq1}\psi^{\prime}(\sigma(\xi_{m}),\sigma
(\xi_{k+m})),\\
\psi_{k}^{\ast}  &  =\sup_{m\geq1}\psi^{\ast}(\sigma(\xi_{m}),\sigma(\xi
_{k+m}))\text{, }\rho_{k}=\sup_{m\geq1}\rho(\sigma(\xi_{m}),\sigma(\xi
_{k+m})).
\end{align*}
By Theorem 7.4 (c and d) in Bradley (2007)%
\begin{align*}
\rho_{k+m}  &  \leq\rho_{k}\rho_{m}\\
1-\psi_{k+m}^{\prime}  &  \leq(1-\psi_{k}^{\prime})(1-\psi_{m}^{\prime}).
\end{align*}
So, by Theorem 7.5(c) in Bradley (2007), if there is $n\geq1$ such that
$\psi_{n}^{\prime}>0,$ then there is $c>0$ such that $1-\psi_{n}^{\prime
}=O(\mathrm{e}^{-cn})\rightarrow0$ as $n\rightarrow\infty.$

Notice that, in terms of conditional probabilities, we also have the following
equivalent definitions:
\begin{align*}
\rho_{1}(\mathbf{\xi)}  &  =\sup_{k}||Q_{k}||_{2}\\
\psi_{1}^{\ast}(\mathbf{\xi)}  &  =\sup_{k}\mathrm{ess}\sup_{x}\sup
_{A\in\mathcal{B}(S)}Q_{k}(x,A)/P_{k}(A)\text{.}\\
\psi_{1}^{\prime}(\mathbf{\xi)}  &  =\sup_{k}\mathrm{ess}\inf_{x}\inf
_{A\in\mathcal{B}(S)}Q_{k}(x,A)/P_{k}(A)\text{.}%
\end{align*}
Note that, by (\ref{cond0})\ we can take $a=\psi_{1}^{\prime}(\mathbf{\xi)}>0$
and $b=\psi_{1}^{\ast}(\mathbf{\xi)}<\infty$. In particular we have
$1-\psi_{k}^{\prime}(\mathbf{\xi)}\leq(1-a)^{k}\rightarrow0$ exponentially fast.

If we consider now measurable functions of a Markov chain $\mathbf{X}%
=(g(\xi_{k}))_{k\geq1},$ by the definition of the mixing coefficients, we
notice that $a<\psi_{1}^{\prime}(\mathbf{X})$ and $\psi_{1}^{\ast}%
(\mathbf{X})<b$.

Actually, the $\psi-$mixing coefficient is defined as
\[
\psi(\mathcal{A},\mathcal{B})=\max[\psi^{\ast}(\mathcal{A},\mathcal{B}%
)-1,1-\psi^{\prime}(\mathcal{A},\mathcal{B})].
\]
For a Markov chain of random variables saying that $\psi_{1}<1$ is equivalent
to $\psi_{1}^{\prime}>0$ and $\psi_{1}^{\ast}<2,$ which implies our condition
(\ref{cond0}).

Assume the variables are centered and have finite second moments. Recall that
$\tau_{n}^{2}=\sum_{j=1}^{n}\mathrm{var}(X_{j}),$ $\sigma_{n}^{2}=E(S_{n}%
^{2})$. From Proposition 13 in Peligrad (2012) we know that for functions of
Markov chains
\[
\frac{1-\rho_{1}}{1+\rho_{1}}\leq\frac{\sigma_{n}^{2}}{\tau_{n}^{2}}\leq
\frac{1+\rho_{1}}{1-\rho_{1}}.
\]
By combining this inequality with Lemma \ref{Bradley2020} we obtain, for
$a>0$
\begin{equation}
\frac{a}{2-a}\leq\frac{\sigma_{n}^{2}}{\tau_{n}^{2}}\leq\frac{2-a}{a}.
\label{varineq}%
\end{equation}
$\ \square$

\subsection{Preliminary general local CLT}

Here we give a general local limit theorem. Its proof is based on the
inversion formula for Fourier transform which is a traditional argument for
this type of behavior. Its statement is practically obtained by arguments in
Section 4 in Hafouta and Kifer, (2016).

\begin{theorem}
\label{GlcltLat}Assume that not all the variables have values in some fixed
lattice. Assume that $b_{n}\rightarrow\infty$ and
\begin{equation}
S_{n}/b_{n}\Rightarrow N(0,1). \label{CLT}%
\end{equation}
In addition, for each $L>0$
\begin{equation}
\ \lim_{T\rightarrow\infty}\text{\ }\lim\sup_{n\rightarrow\infty}%
\int\nolimits_{T\leq|u|\leq Lb_{n}}|E\exp(iu\frac{S_{n}}{b_{n}})|du=0\text{ .}
\label{integral}%
\end{equation}
Then, for any function $h$ on $\mathbb{R}$ which is continuous and with
compact support,
\[
\lim_{n\rightarrow\infty}\sup_{u\in\mathbb{R}}|\sqrt{2\pi}b_{n}Eh(S_{n}%
-u)-\exp(-u^{2}/2b_{n}^{2})\int h(u)\lambda(du)|=0.
\]

\end{theorem}

By decomposing the integral in (\ref{integral}) in two on $\{T\leq
|u|\leq\delta b_{n}\}$ and on $\{\delta b_{n}\leq|u|\leq Lb_{n}\}$ and
changing the variable in the second integral we easily argue that in order to
prove this theorem it is enough to show that for each $L$ fixed there is
$0<\delta<L$ such that%
\[
(C_{1})\text{ \ \ \ }\lim_{T\rightarrow\infty}\lim_{n\rightarrow\infty}%
\sup\int\nolimits_{T\leq|u|\leq b_{n}\delta}|E\exp(iu\frac{S_{n}}{b_{n}%
})|du=0
\]
and
\[
(C_{2})\text{ \ \ \ \ \ \ \ }\lim_{n\rightarrow\infty}b_{n}\int
\nolimits_{\delta<|u|\leq L}|E\exp(iuS_{n})|du=0.\text{ \ \ \ \ \ }%
\]
$\ $

\subsection{Proof of Theorem \ref{ThLocal-lattice}}

For proving Theorem \ref{ThLocal-lattice} we shall verify the conditions of
Theorem \ref{GlcltLat}. The first step is to obtain the CLT. With this aim, we
shall apply Theorem 2.1 in Peligrad (1996). From Bradley (1997), we know that
every lower $\psi-$mixing Markov chain (condition implied by $a>0)$ satisfies
a mixing condition called interlaced $\rho-$mixing, which is precisely the
mixing condition we need to apply Theorem 2.1 in Peligrad (1996). Moreover, by
(\ref{varineq}) and the fact that $a>0,$ condition $\sigma_{n}^{2}%
\rightarrow\infty$ is equivalent to $\tau_{n}^{2}\rightarrow\infty.$ This
means that the Lindeberg's condition (\ref{lindeberg1}) is equivalent to
\[
\frac{1}{\sigma_{n}^{2}}\sum\limits_{k=1}^{n}E(X_{k}^{2}I(|X_{k}%
|\geq\varepsilon\sigma_{n}))\rightarrow0\text{ as }n\rightarrow\infty.
\]
Furthermore, also from (\ref{varineq}), we deduce that
\[
\frac{1}{\sigma_{n}^{2}}\sum\limits_{k=1}^{n}E(X_{k}^{2})\leq\frac{2-a}{a}%
\]
and therefore, all the conditions in Theorem 2.1 Peligrad (1996) are satisfied
and we obtain for this case that
\[
\frac{S_{n}}{\sigma_{n}}\Rightarrow N(0,1).
\]
An alternative way to prove the CLT\ is to use (\ref{varineq}) and then
Theorem 6.48 from Merlev\`{e}de et al. (2019) in the Markov setting.

According to the discussion from the last section, it remains to verify
conditions $(C_{1})$ and $(C_{2})$. We begin by changing the variable in
$(C_{1})$ and, using (\ref{varineq}) and the fact that $a>0,$ we obtain that
$(C_{1})$ is equivalent to
\[
\text{\ \ \ }\lim_{T\rightarrow\infty}\lim\sup_{n\rightarrow\infty}%
\int\nolimits_{T\leq|u|\leq\tau_{n}\delta}|E\exp(iu\frac{S_{n}}{\tau_{n}%
})|du=0.
\]
By Proposition \ref{factor} combined to Condition A, for any $1\leq
|u|<\delta\tau_{n},$

$\ $%
\begin{align*}
|E\exp(iu\frac{S_{n}}{\tau_{n}})|  &  \leq\prod\nolimits_{k=1}^{n}%
[1-\frac{\gamma}{2}(1-|f_{n,k}(\frac{u}{\tau_{n}})|^{2})]^{\frac{1}{4}}\\
&  \leq\exp-\frac{\gamma}{2}\frac{1}{4}\sum\nolimits_{k=1}^{n}(1-|f_{n,k}%
(\frac{u}{\tau_{n}})|^{2})\\
&  \leq\exp(-g(u)).
\end{align*}
Integrating both sides of this inequality on the intervals $T\leq|u|\leq
\delta\tau_{n}$ we obtain
\begin{align*}
\int\nolimits_{T\leq|u|\leq\delta\tau_{n}}|E(iu\frac{S_{n}}{\tau_{n}})|du  &
\leq\int\nolimits_{T\leq|u|\leq\delta\tau_{n}}\exp(-g(u))du\\
&  \leq\int\nolimits_{|u|>T}^{\infty}\exp(-g(u))du.
\end{align*}
Whence, taking first $\lim\sup_{n}$ and then $T\rightarrow\infty,$ condition
$\left(  C_{1}\right)  $ is verified.

We move now to verify $\left(  C_{2}\right)  .$ Because the interval
$[\delta,L]$ is compact, $\left(  C_{2}\right)  $ is verified if we can show
that for any $|u|$ fixed in $[\delta,L]$ we can find an open interval $O_{u}$
such that
\[
\sigma_{n}\sup_{|t|\in O_{u}}|E\exp(itS_{n})|\rightarrow0\text{ as
}n\rightarrow\infty.
\]
By using (\ref{varineq}), it is enough to show that%
\begin{equation}
\tau_{n}\sup_{|t|\in O_{u}}|E\exp(itS_{n})|\rightarrow0\text{ as }%
n\rightarrow\infty. \label{verify}%
\end{equation}
By Proposition \ref{factor}, for any $t$,
\begin{gather*}
\tau_{n}|E\exp(itS_{n})|\leq\tau_{n}\prod\nolimits_{k=1}^{n}[1-\frac{\gamma
}{2}(1-|f_{n,k}(t)|^{2})]^{1/4}\text{ }\\
\leq\tau_{n}\exp(-\ \frac{a^{4}}{8b}\sum\nolimits_{k=1}^{n}(1-|f_{n,k}%
(t)|^{2})\\
\leq\exp[\ln\tau_{n}-\frac{\gamma}{8}\sum\nolimits_{k=1}^{n}(1-|f_{n,k}%
(t)|^{2}).
\end{gather*}
Now (\ref{verify}) is satisfied, provided that
\[
\ln\tau_{n}\left(  1-\inf_{|t|\in O_{u}}\frac{1}{\ln\tau_{n}}\frac{\gamma}%
{8}\sum\nolimits_{k=1}^{n}(1-|f_{k}(t)|^{2})\right)  \rightarrow-\infty.
\]
Since $\tau_{n}\rightarrow\infty$, we obtain in this case that $(C_{2})$
follows from Condition B.

$\ \square$\bigskip

\subsection{Proof of Corollary \ref{corsta-lin}}

First of all we notice that, for any array $\mathbf{X}_{n}=(X_{n,k})_{k}$
defined in this corollary, the mixing coefficients satisfy $a<\psi^{\prime
}(\mathbf{X}_{n})$ and $\psi^{\ast}(\mathbf{X}_{n})<b.$ According to Theorem
\ref{ThLocal-lattice} and the discussion on the conditions A and B, it is
enough to verify Lindeberg's condition\ in (\ref{lindeberg1})\ along to
conditions (\ref{to verify cond A}) and (\ref{to verify B 2}).

Note that $\tau_{n}^{2}=\sum\limits_{k=1}^{n}a_{n,k}^{2}E(X_{0}^{2}%
)\rightarrow\infty$ as $n\rightarrow\infty.$ Hence, by stationarity
Lindeberg's condition becomes, for any $\varepsilon>0,$
\[
\frac{1}{\tau_{n}^{2}}\sum\limits_{k=1}^{n}a_{n,k}^{2}E(X_{0}^{2}I(|X_{0}%
|\geq\varepsilon\tau_{n}))=\frac{E(X_{0}^{2}I(|X_{0}|\geq\varepsilon\tau
_{n}))}{E(X_{0}^{2})}\rightarrow0\text{ as }n\rightarrow\infty,
\]

On the other hand (\ref{to verify cond A})\ becomes: there is $\delta>0$ such
that
\[
\frac{E(X_{0}^{2}I(|a_{n,k}X_{0}|>\delta))}{E(X_{0}^{2})}\leq\frac{E(X_{0}%
^{2}I(|X_{0}|>\delta/M))}{E(X_{0}^{2})}<c\text{ for all }k,
\]
which is obviously satisfied for $\delta$ large enough.

It remains to verify condition (\ref{to verify B 2}). Fix $u\neq0,$ we have to
find an open interval $O_{u}$ containing $|u|$ and a constant $c(u)$ such that
for any $|t|\in U$ we have
\[
|E\exp(ita_{n,k}X_{0})|^{2}=|f_{0}(a_{n,k}t)|^{2}\leq c(u)<1.
\]
As a matter of fact, if $0<c<|u|<d,$ then for any $t$ satisfying $0<c<|t|<d,$
by the boundness of $(a_{n,k}),$ we also have $0<mc<|a_{n,k}t|<Md.$ Now, since
the distribution of $\xi_{0}$ is nonlattice, for any $v$ such that
$0<mc\leq|v|\leq Md$ and because $\ f_{0}$ is continuous on the compact set
$[-Md,-mc]\cup\lbrack mc,Md]$ we can find some constant $C(c,d)$ such that
\[
|f_{0}(v)|\leq C(c,d)<1.
\]
$\square$

\bigskip

\subsection{Proof of Corollary \ref{corlinshort}}

By Theorem 5 in Peligrad and Utev (2006), we know that
\[
\frac{1}{v_{n}}S_{n}\Rightarrow AN(0,1)\text{ as }n\rightarrow\infty,
\]
We also have
\[
\frac{E(S_{n}^{2})}{v_{n}^{2}}\rightarrow A^{2}\text{ as }n\rightarrow\infty,
\]
whence, by (\ref{varineq}), we can find two constants $c_{1}>0$ and $c_{2}>0$
such that \ $c_{1}n\leq E(S_{n}^{2})\leq c_{2}n.$

Recall that $A_{i}=a_{1}+\cdots+a_{i}$ and write
\[
S_{n}=\sum_{i\geq2}b_{n,i}X_{i},
\]
where we used the notation $b_{n,i}=A_{i-1}$ for $2\leq i\leq n$ and
$b_{n,i}=a_{i-n}+\cdots+a_{i-1}=A_{i-1}-A_{i-n-1}$ for $i\geq n+1$.

Let $K_{n}$ be a positive integer such that $K_{n}\geq n$ and $n^{3/2}%
\sum_{\ell\geq K_{n}}|a_{\ell}|\rightarrow0$, as $n\rightarrow\infty$. Let
${\tilde{S}}_{n}=\sum_{i=2}^{n+K_{n}}b_{n,i}X_{i}$. For any $t\in{\mathbb{R}}%
$,%
\begin{gather*}
|E(\exp(itS_{n}))-E(\exp(it{\tilde{S}}_{n}))|\leq2|t|E|S_{n}-{\tilde{S}}%
_{n}|\\
\leq2|t|E|X_{0}|\sum_{i\geq n+K_{n}}|b_{n,i}|\leq2n|t|E|X_{0}|\sum_{\ell\geq
K_{n}}|a_{\ell}|\,.
\end{gather*}
Hence, for each $L>0$ and any sequence $(b_{n})$,
\begin{gather*}
\ \lim_{T\rightarrow\infty}\text{\ }\lim\sup_{n\rightarrow\infty}%
\int\nolimits_{T\leq|u|\leq Lb_{n}}|E\exp(iu\frac{S_{n}}{b_{n}})|du\\
=\ \lim_{T\rightarrow\infty}\text{\ }\lim\sup_{n\rightarrow\infty}%
\int\nolimits_{T\leq|u|\leq Lb_{n}}|E\exp(iu\frac{{\tilde{S}}_{n}}{b_{n}%
})|du\,.
\end{gather*}
Applying now Proposition \ref{factor} to ${\tilde{S}}_{n}=\sum_{i=2}^{n+K_{n}%
}b_{n,i}X_{i}$,
\begin{gather*}
|E(\exp(iu{\tilde{S}}_{n})|^{4}\leq\prod\nolimits_{j=1}^{n+K_{n}}%
[1-\frac{\gamma}{2}(1-|f(b_{n,j}u)|^{2})]\\
\leq\prod\nolimits_{j=1}^{n-1}[1-\frac{\gamma}{2}(1-|f(A_{j}u)|^{2})]\leq
\exp[-\frac{\gamma}{2}\sum\nolimits_{j=1}^{n-1}(1-|f(A_{j}u)|^{2})]\text{ .}%
\end{gather*}
From now on we can proceed exactly as in the proof of Corollary
\ref{corsta-lin}. Indeed, the proof is reduced to verify condition
(\ref{lindeberg1}) and to establish Conditions A and B via the observation
that $m<|A_{k}|<A$ for all $k\geq1$. $\square$

\subsection{Proof of Corollary \ref{corstationary infvar}}

Its proof is based on the next proposition whose prove is similar to that of
Theorem \ref{ThLocal-lattice} and is left to the reader.

In the next proposition $(X_{k})_{k\in\mathbb{Z}}$ is as in Theorem
\ref{ThLocal-lattice}, with the exception that we do not assume that $X_{k}$
has finite second moment. For this case we have:

\begin{proposition}
\label{PrLocal-lattice infvar} Assume that there is a sequence of constants
$b_{n}\rightarrow\infty$ such that
\begin{equation}
\frac{S_{n}}{b_{n}}\Rightarrow N(0,1). \label{conv}%
\end{equation}
Assume that Conditions A and B are satisfied with $\tau_{n}$ replaced by
$b_{n}$ and that$\ $condition (\ref{cond0}) holds. Then, for any function $h$
on $\mathbb{R}$ which is continuous and with compact support,
\begin{equation}
\lim_{n\rightarrow\infty}\sup_{u\in\mathbb{R}}|\sqrt{2\pi}b_{n}Eh(S_{n}%
-u)-\exp(-u^{2}/2b_{n}^{2})\int h(u)\lambda(du)|=0. \label{inf}%
\end{equation}

\end{proposition}

We should notice that Maller (1978), on the pages 106-107, verified Condition
A (with $\tau_{n}$ replaced by $b_{n})$ under the assumptions: for every
$x>0$
\begin{equation}
\sup_{1\leq j\leq n}P(|X_{j}|>b_{n}x)\rightarrow0\text{ as }n\rightarrow
\infty\label{UAN}%
\end{equation}
and

\bigskip

\textbf{Condition \~{A}}$_{1}.$ Denote $V_{n}^{2}(x)=\sum\limits_{k=1}%
^{n}E\left[  (X_{k}-E(X_{k}I(|X_{k}|\leq x))^{2}I(|X_{k}|\leq x)\right]  .$
There are constants $c>0$, $n_{0}\in N$ and $\delta>0$ such that for all
$n>n_{0}$ and $x>\delta$ we have%

\[
\frac{x^{2}\sum\limits_{k=1}^{n}P(|X_{k}|>x)}{V_{n}^{2}(x)}\leq c\text{ }.
\]

\subsection{Proof of Corollary \ref{corstationary infvar}}

We shall verify the conditions in Proposition \ref{PrLocal-lattice infvar}.
First of all, by Lemma \ref{Bradley2020}, we notice that we can apply Theorem
1 in Bradley (1988). Alternatively, one can also use Theorem 2.1 in Peligrad
(1990). It follows that we can find a sequence of positive constants
$b_{n}\rightarrow\infty,$ such that%
\[
\frac{S_{n}}{b_{n}}\Rightarrow N(0,1).
\]

It is well known that saying that $H(x)$ is a slowly varying function as
$x\rightarrow\infty$ is equivalent to
\begin{equation}
\lim_{x\rightarrow\infty}\frac{x^{2}P(|X|>x)}{H(x)}=0. \label{2}%
\end{equation}
Also clearly, since the variables have mean $0$, $\lim_{x\rightarrow\infty
}E(X_{0}I(|X_{0}|<x)=0,$ hence Condition \~{A}$_{1}$ is satisfied. Obviously
condition (\ref{UAN}) is also satisfied and these two properties are precisely
what Maller (1978, pp 107-108), used to show that Condition A is satisfied,
with $\tau_{n}$ replaced by $b_{n}.$

Now by Theorem 18.1.1 in Ibragimov and Linnik (1971), $b_{n}=n^{1/2}%
h(n),$where $h(n)$ as slowly varying at infinity. So
\begin{equation}
\lim_{n\rightarrow\infty}\frac{b_{n}}{b_{n-1}}=1. \label{3}%
\end{equation}
By the same type of arguments used for showing that Condition B$_{2}$ implies
condition B, starting from (\ref{3})\ we show that Condition B is satisfied
with $\tau_{n}$ replaced by $b_{n}$. The proof of this corollary is now
complete. $\square${}

\subsection{On the relation between Conditions A and A$_{1}$}

\begin{lemma}
\label{LemmaA}Let $(X_{j})_{j\geq1}$ satisfying (\ref{conv}) and (\ref{UAN}).
Then Condition A$_{1}$ implies Condition A.
\end{lemma}

\noindent\textbf{Proof.} It is enough to show that Condition A$_{1}$ implies
Condition (1.2) in Maller (1978) and then apply his proof on pages 107-108.
This condition makes used of the symmetrization method. We shall use the
notations: $\tilde{X}_{k}=X_{k}-X_{k}^{\ast}$ with $(X_{k}^{\ast})$ an
independent copy of $(X_{k}).$ We have to verify there are constants $c$,
$n_{0}\in N$ and $\delta>0$ such that for all $n>n_{0}$ and $x>\delta$ we have%

\[
\frac{x^{2}\sum\limits_{k=1}^{n}P(|\tilde{X}_{k}|>x)}{\sum\limits_{k=1}%
^{n}E(\tilde{X}_{k}^{2}I(|\tilde{X}_{k}|\leq x)}\leq c.
\]
By the Markov inequality the expression in the left hand side is dominated by
\[
\frac{\sum\limits_{k=1}^{n}E(\tilde{X}_{k}^{2}I(|\tilde{X}_{k}|>x)}{2\tau
_{n}^{2}-\sum\limits_{k=1}^{n}E(\tilde{X}_{k}^{2}I(|\tilde{X}_{k}|>x)}.
\]
Now, by a desymmetrization argument and monotonicity this quantity is smaller
than%
\[
\frac{8\sum\limits_{k=1}^{n}E(X_{k}^{2}I(|X_{k}|>x/2)}{2\tau_{n}^{2}%
-8\sum\limits_{k=1}^{n}E(X_{k}^{2}I(|X_{k}|>x/2)},
\]
which is uniformly bounded under Condition A.

\section{Acknowledgement}

This paper was partially supported by the NSF grant DMS-1811373. This paper
was developed concomitantly and independently of a new research monograph by
Dolgopyat and Sarig (2020), who treat the local limit theorem for a different
nonstationary situation. The authors would like to thank Richard Bradley for
the statement and proof of Lemma \ref{Bradley2020} and to Jon Aaronson for
useful discussions.

\end{document}